\PassOptionsToPackage{hyphens}{url}
\documentclass[graybox,xcolor=svgnames,usenames,dvipsnames]{svmult}

\usepackage{type1cm}        % activate if the above 3 fonts are
\usepackage{makeidx}         % allows index generation
\usepackage{graphicx}        % standard LaTeX graphics tool
\usepackage{multicol}        % used for the two-column index
\usepackage[bottom]{footmisc}% places footnotes at page bottom
\usepackage{tabularx}

\usepackage{newtxtext}       %
\usepackage[varvw]{newtxmath}       % selects Times Roman as basic font

\usepackage{booktabs}
\usepackage[colorlinks=true, linkcolor=OliveGreen, urlcolor=Fuchsia, citecolor=blue]{hyperref}

\usepackage[T1]{fontenc}

\makeindex             % used for the subject index
\usepackage[dvipsnames]{xcolor}
\usepackage{comment} % provides \excludecomment
\usepackage{soul}

\newif\ifshowcomments
\showcommentstrue  % Trigger

\ifshowcomments
\else
\excludecomment{mm}
\excludecomment{mf}
\excludecomment{me}
\excludecomment{sif}
\fi

\usepackage{mathtools}
\usepackage{stmaryrd}
\newcommand{\F}{\ensuremath{\mathbb{F}}}
\newcommand{\R}{\ensuremath{\mathbb{R}}}
\newcommand{\fl}{\ensuremath{\mathrm{fl}}}
\newcommand{\srp}[1][p]{\ensuremath{\mathrm{SR}_{#1}}}
\newcommand{\srpr}{\srp[p,r]}
\newcommand{\flup}[2][p]{\ensuremath{\llceil #2 \rrceil_{#1}}}
\newcommand{\fldown}[2][p]{\ensuremath{\llfloor #2 \rrfloor_{#1}}}

\newcommand{\expect}{\mathbb{E}}

\usepackage{tikz}
\usetikzlibrary{calc, positioning}
\usetikzlibrary{decorations.pathreplacing, calligraphy, arrows, arrows.meta}
\usepackage{cleveref}

\usepackage[backend=biber,
            giveninits=true,
            style=numeric,     % closest to spmpsci (numeric, bracketed)
            maxnames=2,        % show up to 4 names
            minnames=1,        % then collapse to “et al.”
            sorting=none]{biblatex}

\AtEveryBibitem{\clearfield{month}}
\AtEveryBibitem{\clearfield{issn}}
\AtEveryBibitem{\clearfield{isbn}}
\begin{document}

\title*{What is New in Stochastic Rounding: a Survey on Theory, Hardware, and Applications}

\titlerunning{A Survey on Stochastic Rounding}
% Use \titlerunning{Short Title} for an abbreviated version of
% your contribution title if the original one is too long
\author{El-Mehdi El Arar, Massimiliano Fasi, Silviu-Ioan Filip, and Mantas Mikaitis}

% Use \authorrunning{Short Title} for an abbreviated version of
% your contribution title if the original one is too long
\institute{El-Mehdi El Arar, Sorbonne University, CNRS, LIP6, Paris, France, \email{mehdi.elarar@lip6.fr} \\
  Massimiliano Fasi, University of Leeds, Leeds, UK, \email{m.fasi@leeds.ac.uk}\\
  Silviu-Ioan Filip, Université de Rennes, Inria, CNRS, IRISA, Rennes, France, \email{silviu.filip@inria.fr}\\
Mantas Mikaitis, University of Leeds, Leeds, UK, \email{m.mikaitis@leeds.ac.uk}}
%
% Use the package "url.sty" to avoid-
% problems with special characters
% used in your e-mail or web address
%
\maketitle

% \begin{mf}
%   Currently, we are almost one full page above the limit, so we will need to trim the text quite a bit.

%   \textbf{Update} I've trimmed a few lines, but we are still about half a page over the limit.

%   Thanks - we will do a final language compaction round before submission.
% \end{mf}

\abstract{\emph{Stochastic rounding} (SR) is a probabilistic method used to round numbers to floating-point and fixed-point representations.
  % Whilst round-to-nearest is a default rounding mode on most floating-point units,
  In length $n$ summation, the worst-case error of SR grows as $\sqrt{n}$ with high probability, unlike for standard modes, like round-to-nearest (RN), which grows as $n$.
  For this reason, the former is increasingly employed in large-scale, low-precision computations as an RN alternative.
  Additionally, SR alleviates \emph{stagnation}, whereby relatively small summands are completely rounded off and do not contribute to the sum.
  We provide an update to [Croci et al., \emph{Roy.~Soc.~Open~Sci.}~9.3 (2022), pp. 1--25], a survey which discusses the development and use of SR between 1949 and 2022, citing over 100 references.
  Since then, there has been a surge of new research, and this update covers almost four years of further progress in applying, analysing, and implementing SR.
  Our main focus is \emph{limited-precision stochastic rounding}, a new variant that fixes the precision of the random numbers used.
  We provide insights into industrial and numerical analysis activities surrounding SR, highlighting the next possible steps in making this rounding mode more widely available in hardware.}

\section{What is Stochastic Rounding}

% We will explain what \emph{limited-precision stochastic rounding} is before delving into a formal definition.
% In order to demonstrate, we will be using decimal representation of the reals, before switching to binary floating point in the next section.

Let us start with a simple numerical example using decimal arithmetic. Suppose we want to round the first 20 digits of $\pi=3.14159265358979323846$ to four decimal places and add up the results.
With round-to-nearest (RN), we have
$
3.142+3.142+3.142+\dots
$
where each successive addition contributes an absolute rounding error of approximately $0.00041$ due to rounding up.
With stochastic rounding (SR), however, there is approximately a $59\%$ chance of rounding up, and a $41\%$ chance of rounding down $\pi$ to four significant digits.
These probabilities can be improved by increasing the precision of random numbers when implementing SR (see \cref{sec:error-analysis}).
An ideal scenario would be \emph{exact SR}, which in this case would require using random numbers with $16$ decimal digits of precision to stochastically round the $20$ digits of $\pi$ with exact probabilities.
The computation with SR may be
$
3.142+3.141+3.142+3.141+3.142\dots
$
with absolute rounding errors of approximately $0.00041$ and $-0.00059265$ cancelling out.
\emph{The key difference between RN and SR is that SR rounds to both sides of the exact value with probabilities that are relative to the distances to the surrounding numbers.}
This is what allows it to cancel errors out and reduce the overall summation error.
When the set of rounding errors generated by RN does not have a zero mean they may not cancel out.
We next turn to a formal definition of SR in binary floating-point arithmetic.

\section{Background}
\label{sec:background}

Let $\F \subset \R$ be a normalised floating-point number system with $p$ digits of precision, and let $x \in \R$.
For simplicity of exposition, here we only consider the case in which $x$ lies between the smallest and largest finite value in $\F$.

A deterministic rounding function is a function $\fl \colon \R \to \F$ that uses a fixed rule to map $x \in \R$ to either of the two \emph{rounding candidates},
$\flup{x} = \min \{ y \in \F : y \ge x\}$ or $\fldown{x} = \max \{ y \in \F : y \le x\}$.
An example is RN with ties-to-even (RNE)~\cite[Sec.~4.3.1]{ieee19}, which returns the rounding candidate closest to $x$, using a predetermined rule to break ties if $x$ falls exactly between $\fldown{x}$ and $\flup{x}$.
%Other rounding modes typically available in hardware are the so-called directed rounding functions: round-down (RD), which always returns $\fldown{x}$, round-up (RU), which always returns $\flup{x}$, and round-to-zero (RZ), which returns $\fldown{x}$ if $x \ge 0$ and $\flup{x}$ otherwise.
In contrast, SR may return either rounding candidate with a certain probability.
The precision-$p$ SR is the function $\srp : \R \to \F$ defined by $\srp(x) = x$, if $x \in \F$, and by
\begin{equation}
  \label{eq:srp-def}
  \srp(x) =
  \begin{dcases}
    \flup{x},\quad &\text{with probability } q(x),\\
    \fldown{x}, &\text{with probability } 1- q(x),
  \end{dcases}\qquad
  q(x) = \dfrac{x - \fldown{x}}{\flup{x} - \fldown{x}},
\end{equation}
otherwise.
% An alternative, simpler definition that relies on the concept of unit in the last place (ulp) exists~\cite[Def.~2.1]{effm25}.
The main motivating factor of SR is $\expect(\srp(x)) = x$.
%\begin{equation*}
%\expect(\srp(x)) = x.
%\end{equation*}

\def\ticksep{0.2}
\def\intsep{0.2}
\def\arrowsep{0.2}
\def\rtwosep{0.2}
\def\compssep{0.3}

\def\barheight{0.4}
\def\barsep{0.5}
\def\barlength{2.8}
\def\hdis{1.9}

\begin{figure}[t]
  \centering
  \begin{tikzpicture}[every node/.style={
      minimum width=0pt,
      minimum height=0pt,
      inner sep=0pt},semithick]

    % \node [align=left] at (-0.7,\linesep) {\makebox[0.4cm]{\hfill#8}};

    \node (A) at (0,0) {};
    \node (B) at (3,0) {};
    \node (C) at (4,0) {};
    \node (D) at (8,0) {};

    % Line.
    \draw ($(A)-(\compssep,0)$)--($(D)+(\compssep,0)$);

    % Ticks.
    \foreach \x in {A,B,C,D} % Removed A from this list
    \draw ($(\x)+(0,-\ticksep)$) -- ($(\x)+(0,\ticksep)$);
    % \draw[densely dotted] ($(C)+(0,-\ticksep)$) -- ($(C)+(0,\ticksep)$);

    % Arrows
    % \draw [stealth-stealth,dash pattern={on 0.6pt off 0.6pt}]($(B)+(0,\linesep)$) -- ($(D)+(0,\linesep)$);
    % \draw [->,densely dotted]($(C)+(#6+0.2,\linesep)$) -- ($(D)+(0,\linesep)$);

    % Labels.
    \node [above=\ticksep+0.1 of A] {$\llfloor x \rrfloor_p$};
    \node [above=\ticksep+0.1 of B] {$x$};
    \node [above=\ticksep+0.1 of C] {$\fl_{p+r}(x)$};
    \node [above=\ticksep+0.1 of D] {$\llceil x \rrceil_p$};

    % Intervals.
    % \draw [Stealth-Stealth]($(A)+(0.5*\intsep,-2*\intsep)$) --
    % ($(B)+(-0.5*\intsep,-2*\intsep)$);
    % \node [below=2.6*\intsep of $(A)!0.5!(B)$] {\small$q(x) \cdot  (\flup{x} - \fldown{x})$};
    % \draw [Stealth-Stealth]($(B)+(0.5*\intsep,-2*\intsep)$) --
    % ($(D)+(-0.5*\intsep,-2*\intsep)$);
    % \node [below=2.6*\intsep of $(B)!0.5!(D)$] {\small$\bigl(1-q(x)\bigr) \cdot (\flup{x} - \fldown{x})$};

    % \draw [Stealth-Stealth]($(A)+(0.5*\intsep,-6*\intsep)$) --
    % ($(C)+(-0.5*\intsep,-6*\intsep)$);
    % \node [below=6.6*\intsep of $(A)!0.5!(C)$] {\small$q_r(x) \cdot \fpint{x}$};
    % \draw [Stealth-Stealth]($(C)+(0.5*\intsep,-6*\intsep)$) --
    % ($(D)+(-0.5*\intsep,-6*\intsep)$);
    % \node [below=6.6*\intsep of $(C)!0.5!(D)$] {\small$\bigl(1-q_r(x)\bigr) \cdot (\flup{x} - \fldown{x})$};

  \end{tikzpicture}
  \caption{Difference between the \emph{exact} and \emph{limited-precision stochastic rounding variants}. The former uses probabilities relative to the distances between $x$ and the two precision-$p$ rounding candidates; the latter first rounds $x$ to $\fl_{p+r}(x)$ (where $r$ is the bit-width of the random numbers utilised) and then rounds with probabilities relative to $\fl_{p+r}(x)$. See \eqref{eq:srp-def} and the associated description.}
  \label{fig:sr}
\end{figure}
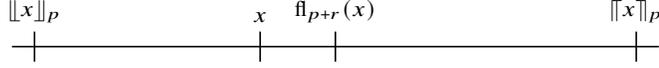

The rounding probability $q(x)$ in~\eqref{eq:srp-def} requires knowledge of the exact $x$, which is typically not available in practical implementations.
This motivates the study of \emph{limited-precision SR}, $\srpr(x) = \srp(\fl_{p+r}(x))$,  where $\fl_{p+r}(x)$ maps $x \in \R$ to precision $p+r$ using some rounding function $\fl_{p+r}$ (\cref{fig:sr}).
%where the rounding probability is computed with respect to $\fl_{p+r}(x)$, an approximation of $x$ with $p + r$ digits.
Since $\fl_{p+r}(x)\neq x$ in general, the probabilities of limited-precision SR diverge from those of exact SR and \mbox{$\expect(\srpr(x)) \neq x$}.

\section{Standardisation of Stochastic Rounding}
\label{sec:standardisation}

The IEEE P3109 interim report~\cite{ieee25} specifies three limited-precision variants of SR
% , each parameterised by an explicit number of random bits
designed to achieve different bias-complexity tradeoffs.
SR only applies within the final projection of the result back to a P3109 format and does not alter the real-arithmetic part of the computation.
It can be used in scalar arithmetic operations, conversions between P3109 formats and between P3109 and IEEE 754 formats, and in block operations.
The exact probability $q(x)$ is approximated by evaluating whether a scaled version of $q(x)$ plus a random integer $0 \le R < 2^{r}$ crosses a threshold.
StochasticA uses the simplest test, comparing $\lfloor 2^{r} q(x) \rfloor + R$ with $2^r$; this reproduces ideal SR in expectation, but the floor operator may introduce a small bias.
StochasticB uses the same test but refines it by using $2^{r+1}$ effective subintervals, which reduces bias but increases complexity of the implementation~\cite{fife25}.
StochasticC replaces the floor operator in StochasticA with RNE, which gives the closest approximation to exact SR~\cite{fife25}.
The document does not specify the number of random bits to be used, or any features of the random number generator.

Tesla’s Configurable Floating Point Formats whitepaper~\cite{tesla} recommends that conversion from bfloat16 and IEEE binary32 to CFloat8 and CFloat16 support both SR and RN.
Arithmetic is typically carried out in higher precision, and SR is confined to the final quantisation step when storing values back into the reduced-precision formats.
% MM: I commented out the below because it requires ML-specific language to be defined, which may be hard in our space limitations.
% The white paper also notes that SR‑based quantisation noise can act as a form of regularisation to prevent overfitting.

\section{Error Analysis of Algorithms with Stochastic Rounding}
\label{sec:error-analysis}

%\cite{ xmh24, xmhk23, effm25}

The numerical analysis of algorithms under SR has attracted increasing attention in recent years for three main reasons. First, probabilistic error models have been developed to better characterise rounding error behaviour. Second, several studies have examined the impact of SR on the stability and accuracy of numerical algorithms in practical applications. Third, the emergence of hardware implementations of SR has motivated analyses that account for implementation constraints and their effect on accuracy.

For a variety of algorithms~\cite{esop22, esop23, esop23a, arar23, oeps25, hi23, des26}, SR yields probabilistic error bounds of order $\mathcal{O}(\sqrt{n}u)$, rather than the classical deterministic bounds of order $\mathcal{O}(nu)$, where $n$ denotes the problem size and $u$ the unit roundoff. To the best of our knowledge, the literature has considered two approaches, both relying on concentration inequalities.

One approach relies on martingale techniques.
First proposed for summation~\cite{chm21}, it has been extended to more complex algorithms and sharper bounds have been obtained through alternative constructions of the martingale. A general methodology, proposed by de Oliveira et al.~\cite{oeps25}, constructs a martingale associated with any direct acyclic graph involving multi-linear errors arising from addition, subtraction, and multiplication.

A more recent approach~\cite{esop23} exploits information on the variance of the error and applies the Bienaymé--Chebyshev inequality to derive tight probabilistic error bounds for a fixed probability and a large $n$. The thesis~\cite{arar23} compares these two approaches~in~detail.

SR is highly effective for machine learning applications (\cref{sec:app}), which are typically based on iterative algorithms such as gradient descent (GD). Xia et al.~\cite{xmhk23} investigated the influence of rounding errors from various sources on the convergence of fixed-stepsize GD for convex problems, showing that SR helps prevent stagnation in low precision. They introduce and analyse two biased variants of SR,
% ($\SR_{\epsilon}$ and signed-$\SR_{\epsilon}$),
which can yield faster GD convergence than standard SR in low precision. These results were subsequently extended to fixed-point arithmetic~\cite{xmh24} for problems satisfying the Polyak--\L{}ojasiewicz inequality. For training Large Language Models (LLMs), Ozkara et al.~\cite{oyp25} presented an error analysis of implicit regularisation and convergence under the Adam optimiser. They showed that Adam’s tolerance can absorb the quantisation noise induced by SR without affecting convergence, particularly when the learning rate is sufficiently large. %Moreover, SR-based updates act as an implicit form of quantisation-aware training, improving robustness to low precision.

%The error analysis was also conducted to support the implementation of SR in hardware.
El Arar et al.~\cite{effm25} developed a new model for the probabilistic error analysis of algorithms using $\srpr$, which they applied to recursive summation and inner product computation of length $n$. Their results suggest that choosing $r$ close to $\lceil (\log_2n) / 2\rceil$ offers the best trade-off between cost and accuracy. This model was further used to analyse Horner's method and pairwise summation, confirming the earlier heuristic~\cite{effm26}.

% \begin{mm}
  % Cite Mehdi's new paper (\url{https://hal.science/% hal-05653338v1}) here and add a short description.
% \end{mm}

%sources: the errors obtained in evaluating the gradient, in computing the multiplication of the rounded gradient with the stepsize, and in determining the subtraction

%\begin{mm}
 % Cite the preprint version of the companion ENUMATH paper on Horner and pair-wise sum err analysis.
%\end{mm}

\section{Software Implementations}

Software-based emulation is useful when hardware SR is not available.
StochasTorch\footnote{\url{https://github.com/nestordemeure/stochastorch}} and
Jochastic\footnote{\url{https://github.com/nestordemeure/jochastic}}
offer support for SR-rounded addition operations in PyTorch and JAX, respectively, with the main use case being parameter updates in neural network training workflows.
Also targeting PyTorch, the mptorch\footnote{\url{ https://github.com/mptorch/mptorch}} package supports limited precision SR, in line with~\cite{effm25}.
For general purpose computations, the Gfloat\footnote{\url{ https://github.com/graphcore-research/gfloat}} Python package operates with generic floating point formats and implements various flavors of limited precision SR operators, with vectorised support for Numpy, JAX, and PyTorch arrays.
srfloat\footnote{\url{ https://github.com/sfilip/srfloat}} supports limited-precision SR addition and multiplication with arbitrary-precision floating-point encodings, offering both Python and C++ interfaces for CPU backends.
LoFloat\footnote{\url{https://github.com/ SudhanvaKulkarni123/LoFloat}} is a C++ library for simulating custom-precision formats and it supports SR.

\section{Hardware}
\label{sec:hardware}

% \section{Hardware in the Research Literature}

%\begin{mm}
%  Target: half a page.
%\end{mm}

\textbf{In the Research Literature.} Chang et al.~\cite{cyls23} presented an SR unit for rounding 16-bit to 8-bit integers with a 3-bit modified Linear-Feedback Shift Register (LFSR) PRNG.
They trained with SR neural networks with 8-bit integer weights, activations, and gradients.
They instantiated 4,096 SR units on an FPGA, finding that the design requires $3.75\times$ fewer LUTs than an 8-bit LFSR, without affecting the neural network's~accuracy.

Ali, Filip, and Sentieys~\cite{afs24} used an LFSR to implement SR in a mixed-precision multiply-accumulate unit with 8-bit floating-point inputs, 12-bit floating-point accumulation, and 18-bit random numbers.
Instead of computing the sum of the significands and then adding random bits to the result~\cite[Sec.~7.3]{cfhm22}, their algorithm adds the random bits as soon as the significand alignment is completed.
This exploits the fact that some of the smaller addend's significand bits will be shifted out to the right and will not be affected by the main addition operation---meaning that the majority of the SR logic can be evaluated in parallel with the shift and add.

Yuan et al.~\cite{ycjl22} explored SR in neural network training using 6-bit arithmetic and 6-bit random numbers.
They argue that one can extract, from within the neural network, random bits that are of higher quality than those obtained from a typical LFSR generator.
This surprising finding eliminates the need for a PRNG in their target machine learning application.
On an FPGA, their design yielded a 38.5\% LUT saving in the SR unit and a 9\% saving overall.
Mishra et al.~\cite{mrc25} explored training of neural networks using an 8-bit LFSR for SR in the conversion operation between \emph{posit} formats.

\medskip

% \section{Commercial Hardware}
% \label{sec:comm-hw}

%\begin{mm}
%  Target: one page.
%\end{mm}
\textbf{Commercial Hardware.} We discuss what is known about the SR-enabled devices announced since the publication of~\cite[Sec.~7.5]{cfhm22}.
% Since the 2022 survey~\cite[Sec.~7.5]{cfhm22}, several new commercial devices with SR have been announced, and here we discuss what is known about them.
\Cref{table:sr-prec} lists the bit-widths of random numbers used in SR by several hardware units for format conversion operations.

Graphcore has released an in-depth description of SR~\cite[Sec.~2.14.6]{graph23}.
When enabled, SR overrides the rounding mode setting and enables a subset of instructions to round stochastically.
First, a binary32~\cite{ieee19} result is produced using RN, and then a pseudo-random number generator (PRNG) generates a bit stream from which between 13 and 24 bits are added to the right-hand side of the significand of the result, before truncating that same region of the significand from the right-hand-side of the result.
The variable length random bit stream is needed for rounding subnormal values in binary16~\cite{ieee19}, based on the number of leading zeros in the destination subnormal.
Any value lower than $2^{-25}$ (half of the smallest subnormal binary16 value) is mapped to zero.
The operations $|a|+|b|$, $a+b$, $a-b$, $a \times b$, and $a \times b + c\times d$ have built-in SR support.
Note that this is not exact SR because of the initial rounding to binary32 within an operation.
A conversion operation from 16- to 8-bit floating-point formats with 4- and 3-bit significands is available. This uses random bit strings up to 11 bits in length.

\begin{table}[t]
  \caption{Number of random bits $r$ used in SR conversion operations between two floating-point formats. The symbol ``--'' indicates that a conversion operation is not specified. The column labels show significand bit-widths of the source $\rightarrow$ destination floating-point format.
    $\dagger$ - when \texttt{cvt.rs} takes four binary32 inputs to round, NVIDIA PTX ISA does not specify how the two halves of a 32-bit random number are utilised to round two separate pairs of binary32 values.}
  \label{table:sr-prec}
  \centering
  \newcolumntype{C}{>{\hspace{3pt}}c<{\hspace{3pt}}}
  \begin{tabularx}{\linewidth}{XCCCCCC@{}}
    \toprule
    Specification & 24 $\rightarrow$ 11 & 24 $\rightarrow$ 8 & 24 $\rightarrow$ 4 & 24 $\rightarrow$ 3 & 11 $\rightarrow$ 4 & 11 $\rightarrow$ 3 \\
    \midrule
    Graphcore ISA~\cite{graph23} & $13$--$24$ & -- & -- & -- & 7--11 & 8--11 \\
    NVIDIA Blackwell~\cite{nvid25} & 13  & $16$ & Up to$16^\dagger$ & Up to$16^\dagger$ & -- & --  \\
    AMD MI300 ISA~\cite{amd24} & --  & -- & $20$ & $21$  & -- & -- \\
    Intel patent~\cite{mmgf25} & $13$ & $16$ & -- & 21 & -- & $8$\\
    \bottomrule
  \end{tabularx}
\end{table}

As for AMD, the MI300 ISA~\cite{amd24} is the first to introduce SR.
It features two instructions: \texttt{CVT\_SR\_FP8\_F32}, which uses 20 random bits~\cite[p.~362]{amd24} to convert from binary32 to fp8-e4m3~\cite{ocp23a}, and \texttt{CV\_SR\_BF8\_F32}, which uses 21~\cite[p.~363]{amd24} to convert from binary32 to fp8-e5m2.
This is exact SR for normalised binary32 values, but limited-precision SR ($r=20$ or $r=21$) for subnormals, for which exact SR would require up to 24 random bits (see Graphcore's approach above).
The instructions take as input the binary32 value to round and an unsigned 32-bit integer used as random bits; as in Graphcore's SR, the random bits are added to the trailing bits of the binary32~value.
The CDNA4 ISA~\cite{amd25}, successor to the MI300 ISA, adds analogous instructions for the 16-bit destination formats fp16 and bfloat16, and for 4-, 6-, and 8-bit destination formats, it includes variants that scale the input by a power of two before applying SR.

Next, we turn to NVIDIA~\cite{nvid25, nvid26}.
The NVIDIA PTX ISA specifies the rounding modifier \texttt{.rs}~\cite[p.~299]{nvid26}, which performs SR by adding random bits to the bits to be rounded off, and using the carry to choose a rounding direction, as done by Graphcore and AMD.
This modifier is available on the NVIDIA B200 and B300 GPUs (\texttt{sm\_100a} and \texttt{sm\_103a} architectures~\cite[p.~302]{nvid26}).
SR is enabled in the \texttt{cvt} PTX instruction when used in conjunction with the rounding modifier, \texttt{cvt.rs}.
Different variants are available for rounding binary32 to binary16, bfloat16, fp8-e5m2, fp8-e4m3, fp6-e2m3, fp6-e3m2, and fp4-e2m1~\cite{ocp23b}.
Some variants take four inputs: the destination register where the results are packed, the two binary32 values to round, and a 32-bit integer containing random bits.
Some of the variants take six inputs instead, with an additional two binary32 inputs.
See \cref{table:sr-prec} for a summary of random number bit-widths used for 8- and 16-bit destination formats~\cite[Sec.~9.7.9.21]{nvid26}.

Luo et al.~\cite{lzwl24} describe Huawei's proposal to round from 24- to 3-bit significands.
They use the 14 least-significant bits of the input significand as random bits, assuming these contain uniformly distributed data.
In the case of fp16 or bfloat16 inputs, $r=2$.
In contrast, Graphcore uses \emph{xoroshiro128+}, whilst the AMD MI300 ISA and NVIDIA PTX ISA specify SR instructions that take random numbers as 32-bit inputs.
Thus, equivalent random number streams could be utilised on AMD and NVIDIA, but reproducibility cannot be assured due to current mismatch in SR precision (\cref{table:sr-prec}).

Finally, Google\footnote{\url{https://cloud.google.com/blog/topics/developers-practitioners/why-stochastic-rounding-is-essential-for-modern-generative-ai}} offers native support for SR in the Ironwood TPU, a domain-specific device for AI training and inference.
Details are not available at the time of writing, but SR is used for training in 4- and 8-bit integer, and 8-bit floating-point formats.

% \section{Hardware Patents}
% \label{sec:patents}

%\begin{mm}
%  Target: half a page.
%\end{mm}

\medskip

\textbf{Patents.} We discuss 7 of the 97 patents filed since February 2022 that mention SR.\footnote{\url{https://patents.google.com/?q=(\%22stochastic+rounding\%22)&before=filing:20260201&after=filing:20220201}}

Three of the relevant patents are from Intel.
Parra Osorio et al.~\cite{eofl23} describe a floating-point conversion instruction with SR that takes the random bits as input.
The examples given are of conversion from binary32 to binary16, using 13 random bits, and from fp16 to fp8-e5m2, using 8.
Mellempudi et al.~\cite{mmgf25} expand on this, adding more example conversion operations.
Fu et al.~\cite{fjhg25} deal with the generation of matrices of random numbers and give a conversion operation with SR as an example application.

From AMD, Shah et al.~\cite{sghm23} also discuss conversion operations with SR.
Their method adds a random value to the trailing significand of the high-precision input, rounds the resulting sum, and finally truncates it to the target format.
A notable feature of this design is the fact that the sum is explicitly rounded before truncation.

Chang and Yu~\cite{chyu24} from NVIDIA consider the use of SR to round in parallel blocks of video pixels or neural network weights, with separate LFSR generators seeded with the positions of the block of pixels or the neural network layer.
This makes SR parallel and reproducible without the need to control the seeds.
In the same spirit, Leshem et al.~\cite{lass25} from Mellanox Technologies describe a technique for performing reproducible SR within network switches.
Without random number generation, random-looking bits are computed from the data itself by processing some of the bits with simple logic gates.
As the randomness is generated from the data itself, the computation is fully reproducible---the SR behaviour depends on the data that arrives in a network switch.
The method is generalised for any input and output floating-point formats and for any number of random bits.
Finally, Oltchik et al.~\cite{orbl25} from Mellanox Technologies describe an idea to generate multiple random seeds from a single seed value and to distribute them across the network to perform SR in a reproducible manner.

NVIDIA, Mellanox Technologies, and Huawei all propose to eliminate the PRNG from SR by drawing random-looking bits from the data itself.

%\begin{mm}
%  A second set of eyes on the patents would be good. Silviu will check the Intel patents. Max will check Mellanox Technologies and NVIDIA.
%\end{mm}

\newcommand{\myparagraph}[1]{\noindent\textbf{#1.}\quad}
\section{Applications of Stochastic Rounding}
\label{sec:app}

\myparagraph{Machine Learning}
% \label{sec:ml}
%\cite{npfl23,oyp25,lagc25,cfbs25,typ25,cxzbrzc25, ztkhyt24,nvid25}
Machine learning applications are one of the main driving factors for the renewed interest in SR in recent years~\cite{cfhm22}.
Indeed, the architectural scaling of LLMs toward the trillion-parameter regime is largely due to the switch to lower precision data types, including $16$-bit and $8$ bit formats, and, more recently, microscaled $4$-bit formats such as MXFP4~\cite{typ25} and NVFP4~\cite{nvid25} in the context of mixed-precision training (MPT)~\cite{mnadegghkv17}.
Since 2022, SR has continued to cement its role as one of the key enablers of MPT: in the backward pass, it is used to compute unbiased estimates of parameter gradients, essential for stable training convergence, and in the parameter update phase, it helps prevent stagnation in summation.

MPT assumes that a high-precision (e.g., binary32 or binary16/bfloat16) copy of the model parameters is kept and updated, whereas lower-precision quantised versions of all other signals are used during forward and backward computations. This quantisation plays a critical role in ensuring that the convergence behaviour of MPT is as close as possible to that of high-precision training. When pretraining LLMs with NVFP4 (groups of $16$ e2m1 values sharing an e4m3 scaling factor), Chmiel~et~al.~\cite[Sec.~3.2]{cfbs25} downcast neural gradients and input activations from bfloat16 using SR, which yields unbiased estimates of the parameter gradients stabilizing convergence. RN has lower quantisation error and is used for the other operations, notably in the forward pass, where SR would inhibit convergence.
NVIDIA's whitepaper on NVFP4 LLM pretraining~\cite{nvid25} similarly argues that SR is beneficial only for quantizing neural~gradients.

The problem of high quantisation error (and consequently high quantisation variance) is also crucial in extremely low-precision training, especially for NVFP4 and MXFP4 (groups of $32$ e2m1 values with an e8m0 scaling factor) encodings. While the scaling factors reduce the dynamic range needed to represent tensor values, outliers can still have a strong impact. In this case, Random Hadamard Transforms (RHTs) are a low-overhead solution that redistributes outliers in an approximate Gaussian distribution, making them easier to represent in narrower formats (see~\cite[Sec.~4.2]{nvid25} and the references therein). This technique has quickly become the norm in many $4$-bit MPT approaches~\cite{nvid25,typ25,cxzbrzc25}. Tseng~et~al.~\cite{typ25} use MXFP4 encodings with RHTs and SR for the backward pass quantisers and note that RHTs act as a variance reducer in SR-quantised GEMMs~\cite[Thm.~3.2]{typ25}. It is similarly used for both backward matrix multiplications by Chen et al.~\cite{cxzbrzc25}, whereas NVIDIA et al.~\cite{nvid25} use it only for the computation of parameter gradients (stating it as optional for the other matrix multiplications).

In microscaling formats, quantisation and scaling operations should be applied along the dot-product dimensions---the inner dimension in a matrix multiplication operation---as this is more hardware friendly. The backward pass requires transposed tensors, thus each tensor must have two different quantised representations for forward and backward pass. This breaks the chain rule, which is believed to reduce model accuracy~\cite[Sec.~4.3]{nvid25} and introduce gradient bias~\cite[Sec.~3.4]{cxzc25}.
To ensure a consistent quantisation in both passes, NVIDIA et al.~\cite{nvid25} propose a 2D $16\times 16$ scaling factor for model parameter tensors in NVFP4 MPT.
Nascimento~et~al.~\cite{npfl23} advocate a similar approach, called direct block floating point, in the context of block floating-point encodings for training convolutional neural networks.
% Compared to~\cite{nvid25},
They use 2D quantisation for input activation and neural gradient tensors, reducing the number of required quantisers, which are performed using SR. Chen~et~al.~(with~\cite{cxzc25} for MXFP$4$ and~\cite{cxzbrzc25} for NVFP$4$) argue for keeping the smaller 1D quantisation dimension, which reduces quantisation error, and instead applying an unbiased gradient \emph{double quantisation} procedure based on SR.

In the parameter update phase, parameter gradients can also be used in low precision.
Converting them with SR removes bias and stabilises the descent procedure, but the benefits diminish once the gradients become too small relative to the quantisation noise/variance (see~\cite[Sec.~4.1]{cfbs25} for a theoretical analysis).
%A Graphcore Research blog post~\cite{fitz26} provides a practical illustration of this effect in a low-precision training loop where the FP6 weights are converted to a wider format for the optimiser computation, updated, and then rounded back to FP6.
%Small optimisation examples are used to motivate SR as a means of preserving small updates in low-precision training.
%However, in a larger example, where the author trains a d12 nanochat with FP6 MLP weights, the low-precision run initially tracks a bfloat16 baseline but later drifts away.
Switching to higher precision---for instance, bfloat16---at the end of training may help close the gap with full high-precision training, as encouraged by NVIDIA~\cite[Appx.~D]{nvid25}. Using a similar SR-quantisation setup,
%~\cite{cfbs25}
Liu~et~al.~\cite[Thm.~1]{lagc25} argue that a large batch size can compensate for reduced precision during backpropagation, mitigating SR-induced variance from per-sample operations.

%Stochastic rounding is one of the main pillars that is making large language model training scale to extremely low-precision datatypes, such as 4-bit MXP4~\cite{typ25} and NVFP4~\cite{nvid25} formats. Indeed, in contexts where network parameter updates $\nabla_{w}\mathcal{L}$ are much smaller than the unit roundoff of the target format, swamping can occur with deterministic rounding updates.

%Its main purpose is to compute unbiased gradient estimates for stable and more accurate model updates.
\medskip
\myparagraph{Neuromorphic Computing}
Neuromorphic computing uses biological neuron-like models: the state is represented by trains of spike times and the learning happens through the plasticity of synaptic connections.
Neuron models and plasticity are approximated by Ordinary Differential Equations (ODEs) solved in fixed-precision arithmetic on digital neuromorphic computers or by analogue circuits.
In 2025, several authors have reported the use of SR in neuromorphic computing.

Atoui et al.~\cite{akbs25} address complex plasticity rules in the context of an analogue neuromorphic computer.
While the machine is analogue, with much of the neuronal simulation performed directly in the circuitry, the plasticity rules are approximated by an embedded digital processor.
The authors utilise a $32$-bit \emph{xorshift} PRNG to round some of the integer multiplication results with SR by comparing the bits to be truncated with a $32$-bit string from the PRNG.

Kim et al.~\cite{ksyl25} implement a spiking neural network with plasticity on an FPGA, using $16$-bit fixed-point arithmetic.
They utilised SR to improve the accuracy of $16$-bit multiplications in the neuron and plasticity models.
As PRNG, they use an $8$-bit LFSR generator, adding random bits to the multiplication result to perform SR.

Urbizagastegui et al.~\cite{usw25} use $8$-bit floating-point arithmetic with SR to simulate synaptic plasticity on a digital computer.
The authors implement $3$-bit SR by adding a $3$-bit random number from a Mersenne twister to the $3$ guard bits in an $8$-bit floating-point adder.
They report that stagnation occurs in their application~\cite[Sec.~3]{usw25} when the magnitudes of the two addends differ so much that the smaller addend's significand is fully shifted out to the right of the $3$ guard bits and lost.
They proposed to use as one of the guard bits a stochastic sticky bit, set with probability that becomes smaller as the difference in the exponents of the two addends increases.
The cost of a stochastic bit compared with a wider PRNG is yet to be explored.

\medskip
\myparagraph{ODEs for Weather, Climate and Epidemiological Simulations}
The impact of SR is particularly evident in weather and climate simulations, where long-term integration and chaotic dynamics show the limitations of RN in \mbox{low precision}.

Klöwer et al.~\cite{kcpp23} analyse how finite-precision arithmetic alters the qualitative behaviour of chaotic, dynamical systems and show that SR preserves correct long-term statistics at reduced precision.
RN creates a finite state space that forces trajectories onto short, periodic orbits.
SR introduces unbiased, random rounding errors that prevent artificial stabilisation and allow trajectories to escape periodic cycles.
Numerical experiments on chaotic simulations show that it is more effective to increase system dimension than precision, and that with SR low-precision simulations remain statistically faithful over long times.
This work relies on \texttt{StochasticRounding.jl}, which uses \emph{xoroshiro128+} to generate the $16$-, $32$-, or $64$-bit random numbers used for rounding.

Kimpson et al.~\cite{kpcp23} look at 100-year long climate simulations under the Simplified Parameterisations PrimitivE Equation DYnamic (SPEEDY) atmospheric model using the ``4$\times$CO2'' experiment, where the level of CO$_{2}$ emissions in a pre-industrial run is abruptly quadrupled.
They show that binary16 arithmetic with RN introduces systematic bias and numerical stagnation, while SR reduces mean bias errors in global temperature and precipitation, with results close to binary32 and binary64 baselines.

Bouali et al.~\cite{beg26} investigated the impact of SR on the numerical solution of ODEs using the explicit midpoint method. Considering two epidemiological models, the SIR and SIRS models, they identified three situations in which RN performs poorly in low precision: near the bifurcation threshold, very small initial infected proportion and SIRS settings where oscillatory behaviour appears. In contrast, SR preserves these qualitative properties. They also proposed a mixed precision approach combined with SR, enabling the use of even lower precision formats for parts of the computation while maintaining the accuracy of the numerical solution.

%\begin{mm}
  % Cite Mehdi's new paper (\url{https://hal.science/hal-05653284v1}) here and add a short description.
% \end{mm}

\medskip
\myparagraph{Scientific Computing}
Creavin~\cite{crea25} considers SR as a mitigation strategy for loss of numerical accuracy in reduced-precision scientific computing.
The work integrates SR, implemented using the rounding logic of \texttt{StochasticRounding.jl}, into the DaCe framework and evaluates it across NPBench kernels and components of the ICON climate model.
Numerical experiments show that SR is unbiased and, compared with RN, can reduce error growth in long accumulation chains by up to three orders of magnitude.
The analysis suggests that SR is most effective when rounding bias dominates truncation error but offers limited benefits in kernels without long summation or with numerically sensitive subroutines.
For NPBench, the author tests various PRNG algorithms, concluding that there is no noticeable effect on the statistical properties~of~SR.

\medskip
\myparagraph{Applied Mathematics}
Dexter et al.~\cite{dbmi25} consider the effect of element-wise SR on the smallest singular value of tall-and-thin matrices rounded to lower precision.
The rounding error is modelled as an independent, zero-mean perturbation with bounded magnitude and nontrivial column-wise variance.
If the matrix is sufficiently tall, this perturbation does not concentrate in low-dimensional subspaces, and the smallest singular value of the rounded matrix is bounded away from zero with high probability, even when the original matrix is rank-deficient.
While RN can preserve or induce rank deficiency, SR implicitly regularises the unrounded matrix, improving its conditioning.

Ma, Yu, and Drineas~\cite{myd26} extend this analysis by showing that the regularising effect of SR is not confined to extremely tall-and-thin matrices and that in fact entire clusters of small singular values at the tail end of the spectrum are lifted when using SR. These results support a broader view of SR as a spectral regulariser.

\medskip
\myparagraph{Program Analysis}
In addition to improving the numerical behaviour of a computation directly, SR can be used in automated program analysis.
PyTracer~\cite{cykg23} is a dynamic profiling tool for Python that detects numerical instabilities by repeatedly executing a program under Monte Carlo arithmetic.
The main novelty of PyTracer is the focus on high-level language interpretability, as it operates at the level of Python functions and data structures rather than individual floating-point operations.

\section{Conclusion}

%\begin{mf}
%  I'm not sure this is the best place for the paragraph above.
%  It seems to break the flow of the discussion.
%  Maybe we could move it to the introduction?
%  \begin{mm}
%    The intro did not seem a good place too. I moved it here as a lead into our summary - how does it read?
%  \end{mm}
%\end{mf}

We have reviewed approximately four years of new research and new hardware and software implementations of SR that followed the survey of Croci~et~al.~\cite{cfhm22}.
The main focus is on \emph{limited-precision stochastic rounding} which, due to limited precision of random numbers used in the implementation, diverges from the exact SR that was the focus of research before and around the time of the survey of Croci~et~al.~\cite{cfhm22}.
New error analysis of SR that accounts for limited-precision randomness has been reviewed (\cref{sec:error-analysis}) and the standardisation of three variants of it in the upcoming IEEE P3109 machine learning floating point standard has been discussed (\cref{sec:standardisation}).
In terms of hardware implementation, we have showed in \cref{table:sr-prec} that implementations in commercial devices and patents (\cref{sec:hardware}) use a varying number of random bits or in some cases do not use a random number generator at all, opting instead to draw random-looking bits from the data.
This lack of consistency across implementations makes it difficult to reason about the numerical behaviour of SR across platforms.
We hope that, by documenting these differences and connecting them with recent error analysis, this survey will help inform future hardware designs and encourage vendors to implement SR more transparently and consistently.
%Understanding how inconsistently it is implemented right now (\cref{table:sr-prec}), should help vendors inform decisions for the future implementations.

%\begin{sif}
%Moved the Drineas and Ipsen paragraph here. It matches better with the forward looking / future work tone of this part.
%\end{sif}

In their SIAM News article of 2024, Drineas and Ipsen~\cite{drip24} interpret SR as a hardware source of random perturbations for algorithmic complexity analysis, in analogy with smoothed analysis. They identify reproducibility, PRNG design, and library integration as the key challenges for its wider adoption. In light of this and the other current trends discussed in the previous sections, future work
on SR may include the following:
\begin{enumerate}
\item Determining the minimum number of random bits needed in SR for a variety of applications (\cref{sec:app}). This should help in providing grounded updates to standards that specify SR operations and encourage alignment of future hardware implementations.
\item Developing techniques for determining undocumented details in some of the hardware implementation of SR (\cref{table:sr-prec}).
\item Studying the effects of SR without a random number generator, such as proposed by patents in \cref{sec:hardware}. SR requires uniform random number generation---studying this aspect would be a first good entry point.
\item Studying the implications of the three SR variants in the IEEE P3109 standard (\cref{sec:standardisation}) to machine learning workloads and scientific computing applications that may utilise future P3109 floating-point arithmetic if it appears on data centre GPUs.
\item Studying the impact of including SR in traditional 32 and 64-bit floating-point units, so that it becomes ubiquitous in floating-point operations such as $\div, +, -, \times$.
\end{enumerate}

%
% \begin{acknowledgement}
% If you want to include acknowledgments of assistance and the like at the end of an individual chapter please use the \verb|acknowledgement| environment -- it will automatically be rendered in line with the preferred layout.
% \end{acknowledgement}

%\section*{Appendix}
%\addcontentsline{toc}{section}{Appendix}
%
%

\printbibliography

%\bibliographystyle{spmpsci}
%\bibliography{references}

\end{document}